\documentclass[12pt]{amsart}
\usepackage[T1]{fontenc}
\usepackage{amsmath,amsthm,amssymb,amsfonts}
\usepackage{hyperref}
\usepackage[capitalise,noabbrev]{cleveref}
\usepackage{newtxtext, newtxmath}
\usepackage{tikz}
\usepackage{graphicx}
\usepackage{float}
\usepackage[utf8]{inputenc}
\usepackage[english]{babel}
\usepackage{csquotes}
\usepackage[url=false, doi=false]{biblatex}
\renewbibmacro{in:}{}
\usepackage[margin=1.25in]{geometry}

\numberwithin{equation}{section}

\theoremstyle{plain}
\newtheorem{theorem}{Theorem}[section]
\newtheorem{lemma}[theorem]{Lemma}
\newtheorem{proposition}[theorem]{Proposition}

\newtheorem{corollary}[theorem]{Corollary}

\newcommand{\R}{\mathbb{R}}
\newcommand{\Rle}{\R_\leq^2}
\newcommand{\Rl}{\R_<^2}

\newcommand{\N}{\mathbb{N}}
\newcommand{\Z}{\mathbb{Z}}

\newcommand{\rk}{\operatorname{rank}}
\newcommand{\vect}{\mathsf{vect}}

\newcommand{\Vect}{\mathsf{Vect}}

\newcommand{\drank}{d_{\rk}}
\newcommand{\ddim}{d_{\operatorname{dim}}}

\newcommand{\Wrk}{W_1^{\rk}}

\newcommand{\isom}{\cong}
\newcommand{\abs}[1]{\left\lvert#1\right\rvert}
\newcommand{\norm}[1]{\left\lVert#1\right\rVert}

\DeclareMathOperator{\kmax}{kmax}
\DeclareMathOperator{\rank}{rank}

\DeclareMathOperator{\id}{Id}

\title[A rank-based distance for interval modules]{A rank-based distance for interval modules and Wasserstein stability of persistence landscapes}

\author{Peter Bubenik}
\email{peter.bubenik@ufl.edu}

\author{Wanchen Zhao}
\email{wzhao1@ufl.edu}

\address{Department of Mathematics, University of Florida, Gainesville, FL 32611, USA}

\begin{document}

\begin{abstract}
    Barcodes and persistence diagrams have a canonical one-parameter family of distances called Wasserstein distances. These distances depend on a choice of distance for interval modules. The usual choices are all Lipschitz equivalent. We observe that they may be written in terms of the dimension function. We instead define a distance for interval modules using the rank function. The rank function contains the information of the dimension function and unlike the dimension function, also contains information on persistence. Using this distance, we show that persistence landscapes are stable with respect to the $1$-Wasserstein distance. This gives us a $1$-Lipschitz embedding of persistence diagrams with the $1$-Wasserstein distance into a Banach space given by an $L^1$ function space, and also provides a lower bound for the $1$-Wasserstein distance. In addition, we give a stability theorem for the mapping from chains to barcodes. This result leads to a heuristic for truncating infinite bars. 
\end{abstract}

\maketitle
 
\section{Introduction}

Homology, computed with coefficients in a field, produces a vector space. 
In applications, one often has a one-parameter family of models, whose homology gives a one-parameter family of vector spaces. A motivation for the development of persistent homology was to understand not just the dimensions of these vector spaces, but the consistencies of these vector spaces across a range of parameter values, which is encoded in the ranks of the linear maps between them~\cite{10.1117/12.57059,MR2701013,elz:tPaS}. 

The family of vector spaces and the linear maps between them is called a persistence module. The rank function is a complete invariant of the persistence module: it may be used to reconstruct the persistence module up to isomorphism. In particular, it determines the dimension function. In contrast, the dimension function provides only upper bounds on the ranks of the maps between parameter values. 

Another complete invariant of persistence modules, which summarizes the rank function, is the barcode~\cite{barcode}. 
The barcode is determined by the decomposition of a persistence module into a direct sum of indecomposable modules called interval modules~\cite{MR1310596,Patel:2018}. 

For quantitative comparisons of persistence modules one desires a distance. 
Barcodes have a canonical (i.e.\ functorial) one-parameter family of distances called Wasserstein distances~\cite{BE2022}.  
These distances depend on a choice of distance between interval modules. 
The first such choice was to use the length of the symmetric difference of the corresponding intervals~\cite{barcode,Collins:2004}. 
Other later choices~\cite{bottleneck} differ from this choice by at most a bounded multiplicative factor. That is, they are Lipschitz equivalent (\cref{lem:dq-lipschitz-equivalent}). 

We observe that this distance for interval modules has an easy expression using the dimension function of the interval module \eqref{eq:ddim}. 
Motivated by the superiority of the rank function from the point of view of persistence, we define an analogous distance for interval modules using the rank function~\eqref{eq:drank}. 
This distance can be expressed in several elementary ways (\cref{lem:drank-alternatives}). 

Our rank-based distance is topologically equivalent to the dimension-based distance, but not Lipschitz equivalent to it (\cref{thm:distance-comparison}). 
However, these distances are locally Lipschitz equivalent for interval modules of positive length (\cref{thm:distance-comparison}). 

Using our rank-based distance for interval modules, we prove that the persistence landscape, as a mapping from barcodes with the $1$-Wasserstein distance to the Banach space $L^1(\N \times \R)$, is $1$-Lipschitz (\cref{thm:Lambda-lipschitz}).
That is, persistence landscapes are stable with respect to their $1$-norm and the $1$-Wasserstein distance.
It follows that the persistence landscape gives a $1$-Lipschitz embedding of persistence diagrams with the $1$-Wasserstein distance into the Banach space $L^1(\N \times \R)$ (\cref{thm:lipschitz-embedding}).

Persistence modules are given by taking the homology of chain complexes of vector spaces and linear maps. 
These chain complexes are composed of persistence modules that are free modules. 
In computational settings, one only uses finitely many values, and it is perhaps reasonable to index persistence modules on a bounded interval. 
Under this hypothesis, we show that the mapping from chain complexes to persistence modules is stable with respect to the $1$-Wasserstein distance using our rank-based distance for interval modules (\cref{thm:stability-chains,thm:stability-j-chains}).
Our statements include an optimal bound for this stability (\cref{thm:stability-chains,thm:stability-j-chains}).
Using this bound, we obtain the following heuristic for truncating infinite bars. 
If the parameter values at which the cells in the complex appear are in the interval $[a,b]$, then the infinite bars should be truncated at $b + (b-a)$.

\subsection*{Related work}

The Wasserstein distances between barcodes and persistence diagrams and their properties have extensively studied~\cite{barcode,bottleneck,Cohen-Steiner:2010aa,BubenikVergili,BE2022,MR4768640,Che:2024aa}. Given a pseudometric on interval modules, there is a functorial construction of the one-parameter family of Wasserstein distances~\cite{BE2022}.
Recently, persistence diagrams and their Wasserstein distances have been studied using the optimal transport of measures~\cite{Divol:2021,MR4496687,Bubenik:2024b,Che:2024}.
Barcodes and persistence diagrams have been shown to be stable with respect the Wasserstein distance~\cite{stability,csem:vineyards,Bubenik:2023a}.

Persistence landscapes are a standard tool for mapping barcodes and persistence diagrams to a Banach space to facilitate topological data analysis~\cite{pl,bubenikDlotko,MR4338670,Betthauser:2022aa,Ayhan:2025,Zhao:2025a}.
Persistence landscapes give a $1$-Lipschitz embedding of persistence diagrams with the $\infty$-Wasserstein distance into the Banach space $L^{\infty}(\N \times \R)$~\cite{pl}.
Our \cref{thm:Lambda-lipschitz,thm:Lambda-bar-lipschitz} are much stronger.
If one uses the dimension-based distance on intervals modules, then the mapping given by persistence landscapes from persistence diagrams with the $p$-Wasserstein distance for $p < \infty$ into $L^q(\N \times \R)$ for $q < \infty$ is not H\"older continuous~\cite{stability}.


Persistent homology has been studied in the setting of chain complexes~\cite{Chacholski:2021,CHACHOLSKI2023101938,Bubenik:2021,MR4878036}. 

\subsection*{Outline of paper}

In \cref{sec:background} we provide background on persistence and distances.
In \cref{sec:metrics-interval-modules} we define our rank-based distance for interval modules and compare it to the previously defined distances by determining the boundaries of its open balls.
In \cref{sec:pl} we prove that the persistence landscape is stable with respect to the $1$-Wasserstein distance when using our rank-based distance for interval modules.
We extend results from finite barcodes and finite persistence diagrams to countable barcodes and countable persistence diagrams.
In \cref{sec:chain-complexes}, we show that the persistent homology of chain complexes of persistence modules is stable with respect to $1$-Wasserstein distance when using our rank-based distance for interval modules if the indexing domain is suitably chosen.

\section{Background} \label{sec:background}

\subsection{Persistence}
\label{sec:persistence}

Let $\R$ denote the set of real numbers together with the usual linear order.
We will often consider $\R$ to be a category whose objects are the real numbers and whose morphisms are given by the linear order.
Let $\mathbf{k}$ be field.
Let $\Vect$ denote the category of $\mathbf{k}$-vector spaces and $\mathbf{k}$-linear maps.
Let $\vect$ denote the full subcategory of finite dimensional vector spaces.
The category $\Vect$ is an abelian category.

A \emph{persistence module} is a functor $M: \R \to \vect$.
For $a \leq b$, we have vector spaces $M_a$ and $M_b$ and the linear map $M_{a \leq b}: M_a \to M_b$.
The \emph{dimension} of $M$ (also called its Hilbert function) is the function $\dim M: \R \to \Z$ given by $\dim M(t) = \dim(M_t)$.
Similarly, the \emph{rank} of $M$ is the function $\rank M: \Rle \to \Z$, where $\Rle = \{(x,y) \in \R^2 \ | \ x \leq y \}$, given by $\rank(M)(a,b) = \rank M_{a \leq b}$.
Let $I$ be an interval in $\R$. 
The corresponding \emph{interval module} $M$ has $M_a = \mathbf{k}$ if $a \in I$ and $M_a = 0$ otherwise, and for $a \leq b$ with $a,b \in I$, $M_{a \leq b}$ is the identity map.
We will denote this interval module by $I$, relying on context to make it clear whether we are referring to an interval or the corresponding interval module.

The structure theorem for persistence modules states any  persistence module $M$ is a direct sum of interval modules~\cite{crawley,MR4143378,Gan:2025}.
That is, $M \isom \bigoplus_{\alpha \in A} I_\alpha$.
Let $\hat{\mathcal{I}}$ denote the set of intervals in $\R$.
The \emph{barcode} of $M$ is the nonnegative function $B(M): \hat{\mathcal{I}} \to \Z$ that assigns each interval its multiplicity 
in the direct sum decomposition $M \isom \bigoplus_{\alpha \in A} I_\alpha$.
We often write barcodes as indexed sets, e.g. $B\left(\bigoplus_{\alpha \in A} I_\alpha\right) = \{I_\alpha\}_{\alpha \in A}$.
Let $\overline{\R}_\leq^2 = \{(x,y) \in [-\infty,\infty] \ | \ x \leq y\}$.
The \emph{persistence diagram} of $M$ is the nonnegative function $D(M): \overline{\R}_\leq^2 \to \Z$ that assigns each pair $(x,y)$ the multiplicity of the intervals
in the direct sum decomposition $M \isom \bigoplus_{\alpha \in A} I_\alpha$
with $\inf I = x$ and $\sup I = y$.
Say that barcodes and persistence diagrams are \emph{finite} if their support is finite.

Let $d$ be a pseudometric on the set of bounded intervals. 
For example, $d(I,J)$ could be given by the length of the symmetric difference $(I \cup J) \setminus (I \cap J)$.
Let $1 \leq p \leq \infty$.
Given two barcodes consisting of bounded intervals, 
$B = \{I_k\}_{k=1}^m$ and $B' = \{J_k\}_{k=1}^n$, 
the \emph{$p$-Wasserstein} distance between them, denoted $W_p(B,B')$ is defined as follows~\cite{barcode,Cohen-Steiner:2010aa}. 
Let $\bar{B} = \{I_k\}_{k=1}^{m+n}$, where $I_{m+1} = \cdots = I_{m+n} = \emptyset$
and $\bar{B}' = \{J_k\}_{k=1}^{m+n}$, where $J_{n+1} = \cdots = J_{m+n} = \emptyset$.
Let $\Sigma_{m+n}$ denote the symmetric group on $\{1,\ldots,m+n\}$.
Then $W_p(B,B') = \min_{\sigma \in \Sigma_{m+n}} \norm{(d(I_k,J_{\sigma(k)}))_{k=1}^{m+n}}_p$.
This distance is a pseudometric on the set of finite barcodes consisting of bounded intervals~\cite{barcode,BE2022}.

Given a bounded interval $I$, the \emph{triangle function} on $I$ is the piecewise linear continuous function $\Lambda(I):\R \to \R$ given by $\Lambda(I)(t) = 0$ if $t \not\in I$, $\Lambda(I)(t) = t - \inf(I)$ if $t \in I$ and $t - \inf I \leq \sup I - t$ and $\Lambda(I)(t) = \sup I - t$ if $t \in I$ and $\sup I - t < t - \inf I$.
Given a barcode $B = \{I_j\}_{j=1}^n$ consisting of bounded intervals, let its \emph{persistence landscape} be the function $\Lambda(B): \N \times \R \to \R$ given by $\Lambda(B)(k,t) = \kmax \{ \Lambda(I_j)(t) \}_{j=1}^n$, where $\kmax$ denotes the $k$th largest of the $n$ numbers and is zero if $k > n$~\cite{pl}.

\subsection{Pseudometrics}

A \emph{pseudometric} on a set $X$ is a nonnegative function $d:X \times X \to \R$ such that for all $x,y,z \in X$, $d(x,x)=0$, $d(x,y)=d(y,x)$, and $d(x,z) \leq d(x,y) + d(y,z)$.
The pair $(X,d)$ is called a \emph{pseudometric space}.
If in addition $d(x,y)=0$ implies that $x=y$ then $d$ is called a \emph{metric} and $(X,d)$ is called a \emph{metric space}.

Let $K \geq 0$. A function $f:(X,d) \to (X',d')$ between pseudometric spaces is said to be \emph{$K$-Lipschitz} if for all $x,y \in X$, $d'(f(x),f(y)) \leq K d(x,y)$. The function $f$ is said to be \emph{Lipschitz} if it is $K$-Lipschitz for some $K$.
Two pseudometrics $d$ and $d'$ on a set $X$ are said to be \emph{Lipschitz equivalent} if both identity map $\id:(X,d) \to (X,d')$ and the identity map $\id:(X,d') \to (X,d)$ are Lipschitz. 
That is, there exists $K,K'$ such that $K d(x,y) \leq d'(x,y) \leq K' d(x,y)$.
The pseudometrics $d$ and $d'$ are said to be \emph{locally Lipschitz equivalent} at $x \in X$ if there exists a neighborhood of $x$ for which they are Lipschitz equivalent.

A pseudometric $d$ on $X$ generates a topology whose base is given by the open balls. 
Lipschitz maps of pseudometric spaces are continuous with respect to the induced topologies.
Two pseudometrics $d$ and $d'$ on $X$ are said to be \emph{topologically equivalent} if they generate the same topology. 
The pseudometrics $d$ and $d'$ are \emph{locally topologically equivalent} at $x \in X$ if the open balls centered at $x$ generate the same filter of neighborhoods. 
That is, each $d$-open ball centered at $x$ contains a $d'$-open ball centered at $x$ and vice versa. 
Two pseudometrics are topologically equivalent if and only if they are locally topologically equivalent at all $x \in X$.
Two pseudometrics are topologically equivalent at $x \in X$ if and only if both identity maps are continuous at $x$ and they are topologically equivalent if and only if both identity maps are continuous.

Let $X$ be a set with an equivalence relation $\sim$ and let $d$ be a pseudometric on $X$.
There is an induced pseudometric $\bar{d}$ on the quotient set $X/\sim$ given by $\bar{d}([x],[y]) = \inf \sum_{i=1}^n d(x_i,y_i)$, where the infimum is taken over all $n \geq 1$ and all $x_1,y_1,\ldots,x_n,y_n \in X$ such that $x_1 \in [x], y_1 \sim x_2, \ldots, y_{n-1} \sim x_n, y_n \in [y]$.

A pseudometric space $(X,d)$ is said to be \emph{complete} if all Cauchy sequences in $X$ converge. 
A \emph{completion} of a pseudometric space $X$ is a complete pseudometric space $\overline{X}$ that contains $X$ has a dense subspace.

Let $(X,d)$ be a pseudometric space. For $x,x' \in X$, define $x \sim x'$ if $d(x,x')=0$.
Then induced pseudometric on the quotient $X/\sim$ is a metric called the \emph{Kolmogorov quotient}.

\section{Metrics for interval modules} \label{sec:metrics-interval-modules}

Let $\mathcal{I}$ denote the set of bounded intervals in $\R$.
Let $\mathcal{I}'$ denote the set of interval modules corresponding to bounded intervals.
There is a canonical bijection from $\mathcal{I}$ to $\mathcal{I}'$ sending an interval $I$ to its corresponding interval module.
In this section we will study pseudometrics on $\mathcal{I}'$.

\subsection{Our rank-based metric for interval modules}
\label{sec:drank}

Let $d_\ell$ denote the pseudometric on $\mathcal{I}$ given by the Lebesgue measure (i.e. length) of the symmetric difference of a pair of intervals.
Let $\ddim$ denote the pseudometric on $\mathcal{I}'$ given by 
\begin{equation} \label{eq:ddim}
    \ddim(I,J) = \int_{\R} \abs{\dim I - \dim J}.
\end{equation}
Under the canonical bijection between $\mathcal{I}$ and $\mathcal{I}'$, $d_\ell = \ddim$.

Now let us compare the dimension function and the rank function. First, for any persistence module $M$ and any $t \in \R$, $\rank M(t,t) = \dim M(t)$. 
So, the rank function contains the information of the dimension function. 
Second, consider the intervals $I = [0,2)$, $J = [0,1)$ and $K=[1,2)$ and their corresponding interval modules.
Then we have a short exact sequence $0 \to K \to  I \to J \to 0$.
This sequence is not split. That is, $I \not\isom J \oplus K$.
We have $\rank I \not= \rank J + \rank K$, but $\dim I = \dim J + \dim K$. 
In general, the dimension function does not detect whether a short exact sequence splits. 
One of the main motivations of studying persistence is to capture information not detected by the dimension function.

Based on these observations, we define an analog of $\ddim$ that uses the rank function instead of the dimension function.
For $I,J \in \mathcal{I}'$, let 
\begin{equation} \label{eq:drank}
    \drank(I,J) = \frac{1}{2} \int_{\Rle} \abs{\rank I - \rank J}.
\end{equation}

\begin{lemma}
    $\drank$ is a pseudometric on $\mathcal{I}'$.
\end{lemma}

\begin{proof}
    Let $I,J,K \in \mathcal{I}$.
    By definition, $\drank(I,I) = 0$ and $\drank(I,J) = \drank(J,I)$.
    Since $\abs{\rank I - \rank K} = \abs{\rank I - \rank J + \rank J - \rank K} \leq \abs{\rank I - \rank J} + \abs{\rank J - \rank K}$, it follows that $\drank(I,K) \leq \drank(I,J) + \drank(J,K)$.
\end{proof}

\begin{figure}
    \centering
    \begin{minipage}[c][1\width]{0.4\textwidth}%
    \begin{tikzpicture}[scale=0.8]
        \draw[thick] (-0.5,-0.5) -- (6.5,6.5);
        \draw[thick] (0,0) -- (0,5) -- (5,5);
        \draw[thick] (2,2) -- (2,6) -- (6,6);
        \filldraw[black] (0,0) circle (2pt) node[anchor=west]{$(a,a)$};
        \filldraw[black] (5,5) circle (2pt) node[anchor=west]{$(b,b)$};
        \filldraw[black] (0,5) circle (2pt) node[anchor=east]{$(a,b)$};
        \filldraw[black] (2,2) circle (2pt) node[anchor=west]{$(c,c)$};
        \filldraw[black] (6,6) circle (2pt) node[anchor=west]{$(d,d)$};
        \filldraw[black] (2,6) circle (2pt) node[anchor=east]{$(c,d)$};
    \end{tikzpicture}
    \end{minipage} \quad \quad
    \begin{minipage}[c][1\width]{0.4\textwidth}%
    \begin{tikzpicture}[scale=0.8]
        \draw[thick] (-0.5,0) -- (6.5,0);
        \draw[thick] (0,0) -- (2.5,2.5) -- (5,0);
        \draw[thick] (2,0) -- (4,2) -- (6,0);
        \filldraw[black] (0,0) circle (2pt) node[anchor=north]{$a$};
        \filldraw[black] (5,0) circle (2pt) node[anchor=north]{$b$};
        \filldraw[black] (2,0) circle (2pt) node[anchor=north]{$c$};
        \filldraw[black] (6,0) circle (2pt) node[anchor=north]{$d$};
    \end{tikzpicture}
    \end{minipage}
    \caption{Consider intervals $I = [a,b)$ and $J=[c,d)$ with $a \leq c \leq b \leq d$. 
    Left: the support of $\abs{\rank I - \rank J}$ is the symmetric difference of two right-angled isosceles triangles.
    Right: the region between the curves $\Lambda(I)$ and $\Lambda(J)$ is also the symmetric difference between two right-angled isosceles triangles.}
    \label{fig:drank}
\end{figure}

Let $I,J \in \mathcal{I}$. Note that $d_\ell(I,J)$ is given by the length of the symmetric difference of the supports of the functions $\dim I$ and $\dim I$.
Analogously, we define $d_a(I,J)$ to be given by half of the area of the symmetric difference of the supports of the functions $\rank I$ and $\rank J$.
Finally, we define a distance using the triangle function.
Let $d_\Lambda(I,J) = \int_{\R} \abs{\Lambda(I) - \Lambda(J)}$.
Let $\hat{d}_a(I,J)$ be given by the area of the symmetric difference of the triangles below the graphs of $\Lambda(I)$ and $\Lambda(J)$.
See \cref{fig:drank}.

\begin{lemma} \label{lem:drank-alternatives}
    Under the canonical bijection between $\mathcal{I}$ and $\mathcal{I}'$, we have $\drank = d_a = \hat{d}_a = d_\Lambda$.
\end{lemma}

\begin{proof}
    Let $I,J$ be nonempty bounded intervals.
    Let $a = \inf I, b= \sup I, c= \inf J$, and $d= \sup J$.
    Assume $a \leq c \leq b \leq d$. See \cref{fig:drank}.
    Then $\drank(I,J) = \frac{1}{2}(\frac{1}{2}(b-a)^2 + \frac{1}{2}(d-c)^2 - 2\frac{1}{2}(b-c)^2) = d_a(I,J)$ and
    $d_\Lambda(I,J) = \frac{1}{2}(b-a)\frac{b-a}{2} + \frac{1}{2}(d-c)\frac{d-c}{2} - 2\frac{1}{2}(b-c)\frac{b-c}{2} = \hat{d}_a(I,J)$.

    If $a \leq c \leq d \leq b$ then
    $\drank(I,J) = \frac{1}{2}(\frac{1}{2}(b-a)^2 - \frac{1}{2}(d-c)^2 = d_a(I,J)$ and
    $d_\Lambda(I,J) = \frac{1}{2}(b-a)\frac{b-a}{2} - \frac{1}{2}(d-c)\frac{d-c}{2} = \hat{d}_a(I,J)$.
    If $a \leq b \leq c \leq d$ then
    $\drank(I,J) = \frac{1}{2}(\frac{1}{2}(b-a)^2 + \frac{1}{2}(d-c)^2 = d_a(I,J)$ and
    $d_\Lambda(I,J) = \frac{1}{2}(b-a)\frac{b-a}{2} + \frac{1}{2}(d-c)\frac{d-c}{2} = \hat{d}_a(I,J)$.

    Also, $\drank(I,\emptyset) = \frac{1}{2} \frac{1}{2} (b-a)^2 = d_a(I,\emptyset)$ and $d_\Lambda(I,\emptyset) = \frac{1}{2}(b-a)\frac{b-a}{2} = \hat{d}_a(I,\emptyset)$.
\end{proof}

Using the canonical bijection between $\mathcal{I}$ and $\mathcal{I}'$, we will take $\drank$ to be a pseudometric on $\mathcal{I}$.


\subsection{Distances induced by norms}

The distance $d_\ell$ is part of a one-parameter family of distances. 
Let $1 \leq q \leq \infty$.
Consider bounded intervals $I$ and $J$.
If $I \neq \emptyset$, let $x = \inf I$ and $y = \sup I$, where $-\infty < x \leq y < \infty$.
Similarly, if $J \neq \emptyset$, let $x' = \inf J$ and $y' = \sup J$.
If $I = J = \emptyset$ then define $d_q(I,J) = 0$. 
If $I \neq \emptyset$ and $J = \emptyset$ then define $d_q(I,J) = \norm{(\frac{y-x}{2},\frac{y-x}{2})}_q$.
Similarly, if $I = \emptyset$ and $J \neq \emptyset$ then define $d_q(I,J) = \norm{(\frac{y'-x'}{2},\frac{y'-x'}{2})}_q$.
If $I$ and $J$ are both nonempty then define
\begin{equation*}
    d_q(I,J) = \min\left( \norm{(x-x',y-y')}_q, \norm{\left(\frac{y-x}{2},\frac{y-x}{2}\right)}_q + \norm{\left(\frac{y'-x'}{2},\frac{y'-x'}{2}\right)}_q \right).
\end{equation*}
Note that $d_1 = d_\ell$.
For each $1 \leq q \leq \infty$, $d_q$ is a  pseudometric on $\mathcal{I}$.

\begin{lemma} \label{lem:dq-lipschitz-equivalent}
    Let $1 \leq q \leq q' \leq \infty$.
    The pseudometrics $d_q$ and $d_{q'}$ on $\mathcal{I}$ are Lipschitz equivalent.
\end{lemma}

\begin{proof}
    For $a,b \in \R$, $\norm{(a,b)}_{q'} \leq \norm{(a,b)}_q \leq 2^{\frac{1}{q} - \frac{1}{q}} \norm{(a,b)}_{q'} \leq 2\norm{(a,b)}_q$.
    Thus, the identity map $\id:(\mathcal{I},d_q) \to (\mathcal{I},d_{q'})$ is $1$-Lipschitz and 
    the identity map $\id:(\mathcal{I},d_{q'}) \to (\mathcal{I},d_{q})$ is $2$-Lipschitz.
\end{proof}

\subsection{Intervals and ordered pairs} \label{sec:intervals-ordered-pairs}

Let $\Rle = \{(x,y) \in \R^2 \ | \ x \leq y\}$.
There is an injective map $\Rle \to \mathcal{I}$ given by mapping $(x,y)$ to the closed interval $[x,y]$.
Using this injection, a pseudometric $d$ on $\mathcal{I}$ restricts to a pseudometric on $\Rle$.

For each of the pseudometrics $d$ on $\mathcal{I}$ defined above, $d(I,J) = 0$ if and only if the interiors of $I$ and $J$ are equal.
Therefore, for the Kolmogorov quotient $(\bar{\mathcal{I}},\bar{d})$, there is a canonical bijection between the set $\bar{\mathcal{I}}$ and the set of bounded open intervals in $\R$, and $\bar{d} = d$ under this bijection.
Note that there is also a canonical bijection between $\bar{\mathcal{I}}$ and $\{\emptyset\} \cup \{ [x,y) \subset \R \ | \ -\infty < x < y < \infty \}$, and $\bar{d} = d$ under this bijection.

Furthermore, the Kolmogorov quotient of $(\Rle,d)$ is the quotient metric space $(\Rle/\Delta,\bar{d})$, where 
$\Delta = \{(x,y) \in \R^2 \ | \ x=y\}$.

For $1 \leq q \leq \infty$, 
the $q$-norm on $\R^2$ induces a metric on $\R^2_{\leq} = \{(x,y) \in \R^2 \ | \ x \leq y \}$.
Then there is an induced metric on the quotient $\R^2_{\leq} / \Delta$.
Denote this metric by $\bar{d}_q$.

There is a bijection between $\bar{\mathcal{I}}$ and $\R^2_{\leq}/ \Delta$ given by mapping the empty interval to $\Delta$ and mapping a nonempty bounded interval $I$ to $[(\inf I, \sup I)]$.
Under this bijection, $d_q = \bar{d}_q$.

Let $\Rl = \{(x,y) \in \R^2 \ | \ x < y\}$.
There is an inclusion $\Rl \to \mathcal{I}$ given by $(x,y) \mapsto [x,y)$.
This inclusion allows us to restrict the pseudometric $\drank$ on $\mathcal{I}$ to a metric on $\Rl$.
If $[x,y) \cap [x',y') = \emptyset$ then 
$\drank((x,y),(x',y')) = \frac{1}{4}(y-x)^2 + \frac{1}{4}(y'-x')^2$,
if $[x,y) \supset [x',y')$ then 
$\drank((x,y),(x',y')) = \frac{1}{4}(y-x)^2 - \frac{1}{4}(y'-x')^2$,
and if $x < x' < y < y'$ then
$\drank((x,y),(x',y')) = \frac{1}{4}(y-x)^2 + \frac{1}{4}(y'-x')^2 - \frac{1}{2}(y-x')^2$.

\subsection{Comparison of rank-based and dimension-based distances}
\label{sec:comparison}

Recall that $\mathcal{I}$ denotes the set of bounded intervals.
Let $\mathcal{I}^+$ denote the set of bounded intervals of positive length.
That is, $\mathcal{I}^+ = \{I \in \mathcal{I} \ | \ I \neq \emptyset \text{ and } \sup I - \inf I > 0\}$.
In this section we will compare the pseudometrics $\drank$ and $\ddim$ on $\mathcal{I}$ and $\mathcal{I}^+$. 
We will prove the following.

\begin{theorem} \label{thm:distance-comparison}
    The pseudometrics $\ddim$ and $\drank$ are locally Lipschitz equivalent on $\mathcal{I}^+$ and topologically equivalent on $\mathcal{I}$.
    They are not Lipschitz equivalent on $\mathcal{I}$.
\end{theorem}


\begin{lemma} \label{lem:not-lipschitz}
    \begin{enumerate}
        \item \label{it:lipschitz-a} 
        There is no $K>0$ such that for all $I,J \in \mathcal{I}$, $\drank(I,J) \leq K \ddim(I,J)$.
        \item \label{it:lipschitz-b} 
        There is no $K>0$ such that for all $I,J \in \mathcal{I}$, $\ddim(I,J) \leq K \drank(I,J)$.
    \end{enumerate}
\end{lemma}

\begin{proof}
    Let $b > 0$. 
    Consider $I = [0,b)$ and $J=\emptyset$.
    Then $\ddim(I,J) = b$ and $\drank(I,J) = \frac{1}{4}b^2$.
    For \eqref{it:lipschitz-a}, $\frac{\drank(I,J)}{\ddim(I,J)} = \frac{b}{4} \to \infty$ as $b \to \infty$.    
    For \eqref{it:lipschitz-b}, $\frac{\ddim(I,J)}{\drank(I,J)} = \frac{4}{b} \to \infty$ as $b \to 0$.
\end{proof}

This result may be restated as follows.

\begin{corollary} \label{cor:not-lipschitz}
    The identity maps 
    $\id:(\mathcal{I},\ddim) \to (\mathcal{I},\drank)$ and
    $\id:(\mathcal{I},\drank) \to (\mathcal{I},\ddim)$ 
    are not Lipschitz.
\end{corollary}

\begin{lemma} \label{lem:continuous}
    Let $I_0 \in \mathcal{I}$ such that $I_0 = \emptyset$ or $\sup I - \inf I = 0$.
    Then the identity maps 
    $\id:(\mathcal{I},\ddim) \to (\mathcal{I},\drank)$ and
    $\id:(\mathcal{I},\drank) \to (\mathcal{I},\ddim)$ 
    are continuous at $I_0$.
\end{lemma}

\begin{proof}
    Let $I \in \mathcal{I}$ such that $\ddim(I,I_0) = r$.
    Then $I \in \mathcal{I}^+$.
    Let $x = \inf I$ and let $y= \sup I$.
    Then $y-x = r$.
    Let $I \in \mathcal{I}$ such that $\drank(I,I_0) = r$.
    Then $I \in \mathcal{I}^+$.
    Let $x = \inf I$ and let $y= \sup I$.
    Then $r = \frac{1}{4}(y-x)^2$ and $y-x = 2\sqrt{r}$.
    Therefore each open ball centered at $I_0$ in one of the metrics contains a ball centered at $I_0$ in the other metric.
\end{proof}


\begin{figure}
    \centering
    \includegraphics[width=0.5\linewidth]{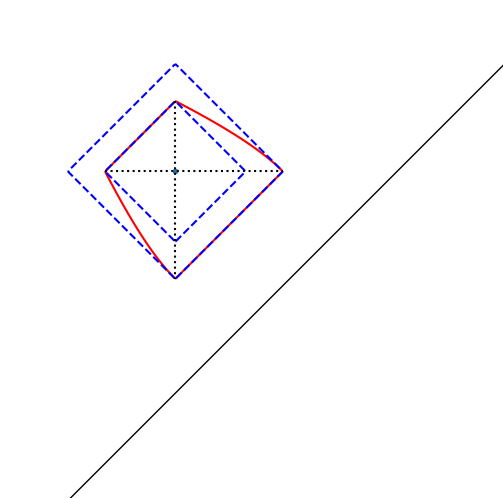}
    \caption{\emph{Nested balls for $\ddim$ and $\drank$.}
    Consider a bounded nonempty interval $I_0$ with positive length and the corresponding point $(\inf I_0, \sup I_0)$.
    We see that a $\ddim$-ball centered at this point with boundary the outer dashed blue curve contains a $\drank$-ball centered at this point with boundary the solid red curve, which in turn contains a $\ddim$-ball centered at this point with boundary the inner dashed blue curve.}
    \label{fig:nested-balls}
\end{figure}

Let $I_0 \in \mathcal{I}^+$.
Let $x_0 = \inf I_0$ and $y_0 = \sup I_0$.
Then $\ell_0 = y_0 - x_0 > 0$.
Furthermore $\ddim(I_0,\emptyset) = \ell_0$ and $\drank(I_0,\emptyset) = \frac{1}{4}\ell_0^2$.
Let $0 < r < \frac{1}{4}\ell_0^2$.
Let $\rho = \frac{4r}{\ell_0^2}$.
Then $0 < \rho < 1$.
Let $s = \ell_0(\sqrt{1+\rho} - 1)$ and $t = \ell_0(1 - \sqrt{1-\rho})$.
For a pseudometric $d$ and $a > 0$, let $B_a^d(I_0)$ denote the $d$-open ball of radius $a$ centered at $I_0$.
See \cref{fig:nested-balls}. 

\begin{lemma}
    $B_s^{\ddim}(I_0) \subset B_r^{\drank}(I_0) \subset B_t^{\ddim}(I_0)$.
\end{lemma}

\begin{proof}
     The result will follow from an analysis of the set of intervals on the boundary of the open ball of radius $r$ centered at $I_0$ for $\ddim$ and $\drank$.


    Let $I \in \mathcal{I}$ such that $\ddim(I,I_0) = r$.
    Then $I \in \mathcal{I}^+$.
    Let $x = \inf I$ and $y = \sup I$.
    There are four cases.
    If $x \leq x_0 < y_0 \leq y$ then $y = x + \ell_0 + r$.
    If $x_0 \leq x < y \leq y$ then $y = x + \ell_0 - r$.
    If $x > x_0$ and $y > y_0$ then $y = - x + x_0+y_0 + r$.
    If $x < x_0$ and $y < y_0$ then $y = -x + x_0+y_0 - r$.
    So the boundary of the $\ddim$ ball is given by a piecewise linear curve consisting of lines of slope $\pm 1$.

    Let $I \in \mathcal{I}$ such that $\drank(I,I_0) = r$.
    Then $I \in \mathcal{I}^+$.
    Let $x = \inf I$, $y = \sup I$, and let $\ell= y-x$.
    There are four cases.
    If $x \leq x_0 < y_0 \leq y$ then $r = \frac{1}{4}\ell^2 - \frac{1}{4}\ell_0^2$.
    From this we get that $y = x + \sqrt{l_0^2 + 4r}$.
    If $x_0 \leq x < y \leq y$ then $r = \frac{1}{4}\ell_0^2 - \frac{1}{4}\ell^2$.
    From this we get that $y = x + \sqrt{l_0^2 - 4r}$.
    If $x > x_0$ and $y > y_0$ then 
    $r = \frac{1}{4}\ell_0^2 + \frac{1}{4}\ell^2 - \frac{1}{2}(y_0-x)^2$.
    From this we get that $y = x + \sqrt{2(y_0-x)^2 - (l_0^2 - 4r)}$.
    If $x < x_0$ and $y < y_0$ then 
    $r = \frac{1}{4}\ell^2 + \frac{1}{4}\ell_0^2 - \frac{1}{2}(y-x_0)^2$.
    So the boundary of the $\drank$ ball is given by a curve consisting of two lines of slope $1$ and the graphs of two quadratic functions.

    The boundaries between these these four cases are given by the following four ``corner'' cases.
    If $x = x_0 < y_0 \leq y$ then $y = y_0 +s$, where $s = \sqrt{\ell_0^2 + 4r} - \ell_0$.
    If $x \leq x_0 \leq y_0 = y$ then $x = x_0 - s$.
    If $x_0 = x < y \leq y_0$ then $y = y_0 - t$, where $t = \ell_0 - \sqrt{\ell_0^2 - 4r}$.
    If $x_0 \leq x < y = y_0$ then $x = x_0 + t$.

    To complete the proof it remains to show that for 
    $x \geq x_0$ and $y \geq y_0$ the graph of
    $y = x + \sqrt{2(y_0-x)^2 - (l_0^2 - 4r)}$ lies between the graphs of 
    $y = -x + x_0 + y_0 + s$ and 
    $y = -x + x_0 + y_0 + t$.
    The case $x \leq x_0$ and $y \leq y_0$ will follow by symmetry.

    Consider the difference $y = x + \sqrt{2(y_0-x)^2 - (\ell_0^2 - 4r)} - \left(-x +2x_0 + \sqrt{\ell_0^2 + 4r}\right)$.
    Note that $y_0 + s = x_0 + \sqrt{\ell_0^2 + 4r}$ and $x_0 + t = y_0 - \sqrt{\ell_0^2 -4r}$.
    When $x = x_0$, $y = x_0 + \sqrt{\ell_0^2 + 4r} - (x_0 + \sqrt{\ell_0^2 + 4r} = 0$.
    Now write the difference as $y = x + \sqrt{2(y_0-x)^2 - (\ell_0^2 - 4r)} - (-x + x_0 + y_0 + s)$.
    When $x = x_0 + t$ then $x = y_0 - \sqrt{\ell_0^2 - 4r}$.
    Thus $y = x_0 + t + \sqrt{\ell_0^2 - 4r} - (-t + y_0 + s) = t-s$.
    We claim that $t-s>0$. 
    Indeed, let $\rho = \frac{4r}{\ell_0^2}$.
    Then $0 < \rho < 1$ and $s = \ell_0(\sqrt{1+\rho}-1)$ and $t = \ell_0(1 - \sqrt{1-\rho})$.
    Hence $t-s = \ell_0(2 - \sqrt{1+\rho} - \sqrt{1-\rho}) > 0$.
    It remains to show that the difference is increasing between $x=x_0$ and $x = x_0 + t$.

    For the difference, we compute the derivative.
    \begin{equation*}
        y' = 2 - \frac{2(y_0-x)}{\sqrt{2(y_0-x)^2 - (\ell_0^2 - 4r)}}
    \end{equation*}
    When $x=x_0$, $y' = 2 - \frac{2\ell_0}{\sqrt{\ell_0^2+4r}} > 0$.

    Finally, we compute the second derivative.
    \begin{equation} \label{eq:y''}
        y'' = \frac{-2(\ell_0^2 - 4r)}{[2(y_0-x)^2 - (\ell_0^2 - 4r)]^{\frac{3}{2}}}
    \end{equation}
    Since $x_0 \leq x \leq x_0 + t$, where $t = \ell_0 - \sqrt{\ell_0^2-4r}$, we have 
    $y_0 - x_0 - t \leq y_0 - x \leq y_0 - x_0$.
    Thus, 
    $\sqrt{\ell_0^2-4r} \leq y_0 - x \leq \ell_0$.
    From this, we get that the denominator of \eqref{eq:y''} lies between 
    $(\ell_0^2-4r)^{\frac{3}{2}}$ and 
    $(\ell_0^2+4r)^{\frac{3}{2}}$,
    and is therefore always positive.
    Hence $y''<0$.
    Thus the sign of $y'$ does not change, and therefore $y$ is increasing as desired.
\end{proof}

\begin{proposition} \label{prop:locally-lipschitz-equivalent}
    The pseudometrics $\ddim$ and $\drank$ are locally Lipschitz equivalent on $\mathcal{I}^+$.
\end{proposition}
    
\begin{proof}
    As $r \to 0$, $\rho \to 0$, 
    $s \to \ell_0 \frac{\rho}{2} =  \frac{2}{\ell_0}r$ and 
    $t \to \ell_0 \frac{\rho}{2} =  \frac{2}{\ell_0}r$.
%
    Hence, for $r$ sufficiently small, 
    $\drank(I,I_0) \leq r$ implies that $\ddim(I,I_0) \leq \frac{4}{\ell_0} r$
    and
    $\ddim(I,I_0) \leq r$ implies that $\drank(I,I_0) \leq \ell_0 r$.

    Let $I_1$ and $I_2$ be sufficiently close to $I_0$. 
    Then the lengths of $I_1$ and $I_2$, $\ell_1$ and $\ell_2$ respectively, satisfy $\frac{1}{2}\ell_0 \leq \ell_1, \ell_2 \leq 2 \ell_0$.
    Let $r = \drank(I_1,I_2)$.
    Then $\ddim(I_1,I_2) \leq \frac{4}{\ell_1}r \leq \frac{8}{\ell_0} \drank(I_1,I_2)$.
    Let $r' = \ddim(I_1,I_2)$.
    Then $\drank(I_1,I_2) \leq \ell_1 r' \leq 2 {\ell_0} \ddim(I_1,I_2)$.
\end{proof}

\cref{thm:distance-comparison} follows from \cref{prop:locally-lipschitz-equivalent}, \cref{lem:continuous}, and \cref{cor:not-lipschitz}.

\subsection{The shape of balls for our rank-based distance}

In this section, we discuss the shape of $\drank$ balls.

The boundaries of balls centered at intervals in $\mathcal{I} \setminus \mathcal{I}^+$ for both $\ddim$ and $\drank$ are given by lines $y = x + c$.
The boundaries of balls centered at intervals in $\mathcal{I}^+$ are more interesting. 

\begin{figure}
    \centering
    \includegraphics[width=0.5\linewidth]{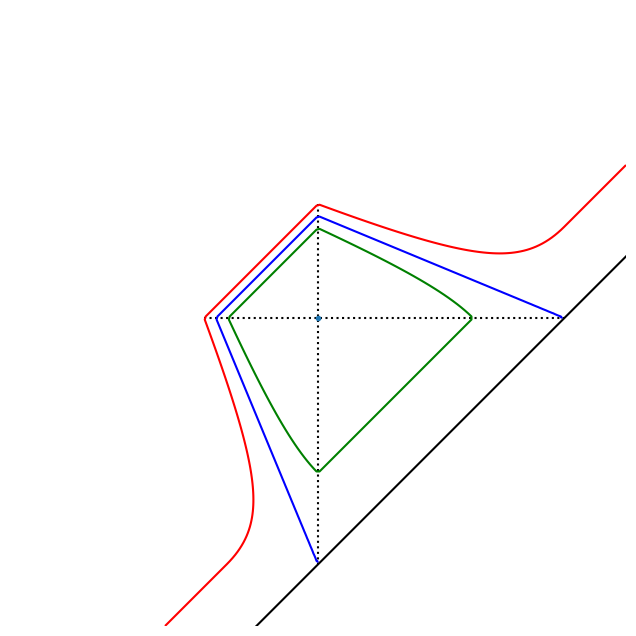}
    \caption{\emph{Balls centered at an interval of positive length.}
    Consider a bounded interval $I$ of positive length and the corresponding point $(\inf I, \sup I)$. 
    There are three cases for the shape of a $\drank$ ball centered at $I$ with radius $r$:
    $r < \drank(I,\emptyset)$; $r = \drank(I,\emptyset)$; and $r > \drank(I,\emptyset)$.
    We plot the boundaries of three balls centered at $I$ for these three cases.
    }
    \label{fig:balls-of-different-shape}
\end{figure}

Let $I_0 \in \mathcal{I}^+$ and let $x_0 = \inf I_0$ and $y_0 = \sup I_0$.
Let $r>0$ and consider the set of all bounded intervals with $\drank(I,I_0) = r$.
There are three cases, depending on whether $r$ is less than, equal to, or greater than $\drank(I_0,\emptyset)$. See \cref{fig:balls-of-different-shape}.
Note that $\drank(I_0,\emptyset) = \frac{1}{4}\ell_0^2$, where $\ell_0 = y_0 - x_0$.
The first case, $r < \frac{1}{4}\ell_0^2$ was considered in \cref{sec:comparison}.
We will consider the remaining two cases.

In the second case, $r = \frac{1}{4}\ell_0^2$.
Whenever $I$ is nonempty, let $x = \inf I$, $y= \sup I$, and $\ell = y-x$.
There are four sub-cases. 
In the first, $I \supset I_0$. 
Then $4r = \ell^2 - \ell_0^2$.
Therefore $y= x + 2 \ell_0^2$.
In the second, $I = \mathcal{I} \setminus \mathcal{I}^+$. 
Therefore the boundary of the circle is given by $y=x$.
In the third, $x_0 < x$ and $y_0 < y$.
Then $4r = \ell^2 + \ell_0^2 - 2(y_0-x)^2$.
Therefore $y = -(\sqrt{2}-1) x + \sqrt{2} y_0$.
In the fourth, $x<x_0$ and $y< y_0$.
Then $4r = \ell^2 + \ell_0^2 - 2(y-x_0)^2$.
Therefore $y = -\frac{1}{\sqrt{2}-1} x + \frac{x_0}{\sqrt{2}-1}$.
Notice that the boundary of the ball is given by a trapezoid with corners
$(x_0,y_0+s)$, 
$(x_0-s,y_0)$,
$(x_0,x_0)$, and
$(y_0,y_0)$,
where $x = (\sqrt{2}-1) \ell_0$.

In the third case, $r > \frac{1}{4} \ell_0^2$.
Whenever $I$ is nonempty, let $x = \inf I$, $y= \sup I$, and $\ell = y-x$.
There are four sub-cases. 
In the first, $I \supset I_0$. 
Then $4r = \ell^2 - \ell_0^2$.
Therefore $y= x + \sqrt{4r + \ell_0^2}$.
In the second, $I \cap I_0 = \emptyset$. 
Then $4r = \ell_0^2 + \ell^2$.
Therefore $y = x + \sqrt{4r - \ell_0^2}$.
In the third, $x_0 < x$ and $y_0 < y$.
Then $4r = \ell^2 + \ell_0^2 - 2(y_0-x)^2$.
Therefore $y = x + \sqrt{2(y_0 - x)^2 + 4r - \ell_0^2}$.
In the fourth, $x<x_0$ and $y< y_0$.
Then $4r = \ell^2 + \ell_0^2 - 2(y-x_0)^2$.
Therefore $x = y - \sqrt{2(y-x_0)^2 + 4r - \ell_0^2}$.

\begin{figure}
    \centering
    \includegraphics[width=0.9\linewidth]{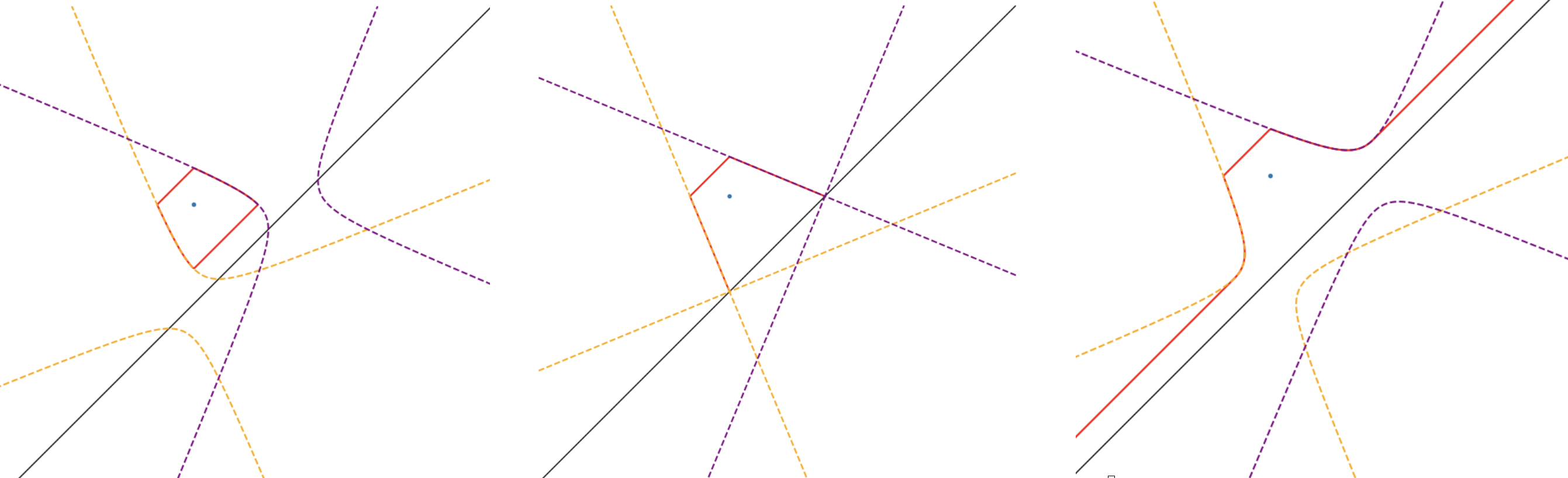}
    \caption{\emph{Boundaries of balls centered at intervals of positive length and their corresponding hyperbolas.}
    We plot the boundaries of the balls in the three cases in \cref{fig:balls-of-different-shape} together with their corresponding hyperbolas.
    }
    \label{fig:hyperbolas}
\end{figure}

See \cref{fig:hyperbolas} for the hyperbolas containing the portion of the boundaries of the balls given by the last two sub-cases.

\section{Stability of persistence landscapes} \label{sec:pl}

Let $\Wrk$ denote the $1$-Wasserstein distance (\cref{sec:persistence}) based on $\drank$ (\cref{sec:drank}).

\subsection{Finite barcodes}

Recall that the persistence landscape is a map from finite barcodes consisting of bounded intervals to real-valued functions on $\N \times \R$. We will show that this map is $1$-Lipschitz with respect to $\Wrk$ and the $1$-norm.

Let $B$ and $B'$ be finite barcodes consisting of bounded intervals. 
\begin{theorem} \label{thm:stability-pl}
    $\norm{\Lambda(B) - \Lambda(B')}_1 \leq \Wrk(B,B')$.
\end{theorem}

\begin{proof}
Let $B  = \{I_\ell\}_{\ell=1}^m$ and $B' = \{J_\ell\}_{\ell=1}^n$ be barcodes consisting of bounded intervals. 
Let $\bar{B}  = \{I_\ell\}_{\ell=1}^{m+n}$ and $\bar{B}' = \{J_\ell\}_{\ell=1}^{m+n}$ where $I_{m+1} = \cdots = I_{m+n} = J_{n+1} = \cdots = J_{m+n} = \emptyset$.
Choose $\sigma \in \Sigma_{m+n}$ such that $\Wrk(B,B') = \norm{(\drank(I_\ell,J_{\sigma(\ell)}))_{\ell=1}^{m+n}}_p$.

Fix $t \in \R$.
For $1 \leq \ell \leq m+n$, let $a_\ell = \Lambda(I_\ell)(t)$ and $b_\ell = \Lambda(J_\ell)(t)$.
For $1 \leq k \leq m+n$, let $a_{(k)} = \kmax(a_1,\ldots,a_{m+n})$ and $b_{(k)} = \kmax(b_1,\ldots,b_{m+n})$.
It is an elementary lemma 
that the optimal way to transport unit masses at $r$ locations in $\R$ to another $r$ locations is given by ordering the locations and moving the rightmost mass to the rightmost location and so on. Therefore,
\begin{equation} \label{eq:OT}
    \sum_{k=1}^{m+n} \abs{a_{(k)} - b_{(k)}} \leq \sum_{\ell=1}^{m+n} \abs{a_\ell - b_{\sigma(\ell)}}.
\end{equation}
The left hand side of \eqref{eq:OT} equals $\sum_{k=1}^{m+n} \abs{\Lambda(B)(k,t) - \Lambda(B')(k,t)}$.
Integrating, we obtain $\norm{\Lambda(B)-\Lambda(B)}_1$.
The right hand side of \eqref{eq:OT} equals $\sum_{\ell=1}^{m+n} \abs{\Lambda(I_\ell)(t) - \Lambda(J_{\sigma(\ell)})(t)}$.
Integrating, we obtain $\sum_{\ell=1}^{m+n} d_\Lambda(I_\ell,J_{\sigma(\ell)}) = \Wrk(B,B')$.
\end{proof}

Recall that $\mathcal{I}$ is the set of bounded intervals in $\R$. 
Let $D(\mathcal{I})$ denote the set of nonnegative functions $\alpha: \mathcal{I} \to \Z$ with finite support. That is, it is the set of formal sums on $\mathcal{I}$. 
We call such formal sums of intervals \emph{finite barcodes}. 
Then $(D(\mathcal{I}),\Wrk)$ is a pseudometric space. 
We now restate \cref{thm:stability-pl}. 

\begin{theorem} \label{thm:Lambda-lipschitz}
    $\Lambda: (D(\mathcal{I}),\Wrk) \to L^1(\N \times \R)$ is $1$-Lipschitz.
\end{theorem}

\subsection{Finite persistence diagrams}

Recall that $\Rle = \{ (x,y) \in \R^2 \ | \ x \leq y \}$.
We have the quotient set $\Rle/\Delta$.
Let $\Rl = \{ (x,y) \in \R^2 \ | \ x < y\}$.
Then $\Rle/\Delta = \Rl \cup \{\Delta\}$.
Let $D(\Rl)$ denote the set of nonnegative functions $\alpha: \Rl \to \Z$ with finite support. 
We call such formal sums \emph{finite persistence diagrams}. 

The inclusion $\Rl \to \mathcal{I}$ given by $(x,y) \mapsto [x,y)$ allows us to restrict the pseudometric $\drank$ on $\mathcal{I}$  to a metric on $\Rl$.
In addition, this inclusion induces an inclusion $D(\Rl) \to D(\mathcal{I})$.
This inclusion allows us to restrict the pseudometric $\Wrk$ to a metric on $D(\Rl)$.

\begin{theorem} \label{thm:lipschitz-embedding}
    $\Lambda: (D(\Rl),\Wrk) \to L^1(\N \times \R)$ is a $1$-Lipschitz embedding.
\end{theorem} 

\begin{proof}   
    The map $\Lambda: (D(\Rl),\Wrk) \to L^1(\N \times \R)$ is $1$-Lipschitz by \cref{thm:Lambda-lipschitz}.
    Furthermore, the persistence landscape is injective for finite persistence diagrams~\cite{pl}.
\end{proof}
That is, the persistence landscape is a $1$-Lipschitz embedding of finite persistence diagrams on $\Rl$ into the Banach space $L^1(\N \times \R)$.

\subsection{Countable barcodes}

Let $\overline{D}(\mathcal{I})$ be the set of nonnegative functions $\alpha: \mathcal{I} \to \Z$ such that 
$\sum_{I \in \mathcal{I}} \alpha(I) \drank(I,\emptyset) < \infty$.
We call such functions \emph{countable barcodes}.
It follows from this condition that $\alpha = \sum_{j \in J} I_j$ for some countable set $J$.
Furthermore, if $\alpha = \sum_{j=1}^\infty I_j$ and
$\sum_{j=1}^\infty \drank(I_j,\emptyset) < \infty$,
then 
it is easy to check that
$(\sum_{j=1}^n I_j)$ is a Cauchy sequence in $(D(\mathcal{I}),\Wrk)$.

We extend $\drank$ from $D(\mathcal{I})$ to $\overline{D}(\mathcal{I})$ as follows.
For $\alpha = \sum_{j=1}^\infty I_j, \beta = \sum_{j=1}^\infty I'_j \in \overline{D}(\mathcal{I})$,
define $\Wrk(\alpha,\beta) = \lim_{n \to \infty} \Wrk( \sum_{j=1}^n I_j, \sum_{j=1}^n I'_j)$.

\begin{lemma}[{\cite[Theorem 6.20, Lemma 6.17]{MR4496687}}]
    $\Wrk$ is well defined on $\overline{D}(\mathcal{I})$ and $(\overline{D}(\mathcal{I}),\Wrk)$ is a completion of $(D(\mathcal{I}),\Wrk)$.
\end{lemma}

\begin{theorem} \label{thm:Lambda-bar-lipschitz}
    $\Lambda$ extends to a $1$-Lipschitz map $\overline{\Lambda}: (\overline{D}(\mathcal{I}),\Wrk) \to L^1(\N \times \R)$.
\end{theorem}

\cref{thm:Lambda-bar-lipschitz} will follow immediately from \cref{thm:Lambda-lipschitz} and the lemma below.

Let $f:X \to Y$ be a $K$-Lipschitz map from a pseudometric space to a complete metric space.
Let $\overline{X}$ be a completion of $X$.
Define an extension $\overline{f}$ of $f$ to $\overline{X}$ as follows.
For $x \in \overline{X}$, choose a Cauchy sequence $(x_n)$ in $X$ that converges to $x$.
Let $\overline{f}(x)$ be given by the limit of the convergent sequence $(f(x_n))$.

\begin{lemma}
    The function $\overline{f}: \overline{X} \to Y$ is $K$-Lipschitz.
\end{lemma}
    
\begin{proof}
    Let $x,y \in \overline{X}$.
    Choose Cauchy sequences $(x_n)$ and $(y_n)$ in X that converge to $x$ and $y$, respectively.
    Then $d_Y(f(x),f(y)) = d_Y(\lim f(x_n), \lim f(y_n)) = \lim d_Y(f(x_n),f(y_n)) \leq \lim K d_X(x_n,y_n) = K d_X(x,y)$.
\end{proof}

\subsection{Countable persistence diagrams}

We define $\overline{D}(\Rl)$ to be the set of nonnegative functions $\alpha: \Rl \to \Z$ such that $\sum_{(x,y) \in \Rl} \frac{1}{4}(y-x)^2 < \infty$. 
As above, we extend $\Wrk$ from $D(\Rl)$ to $\overline{D}(\Rl)$.
Then $(\overline{D}(\Rl),\Wrk)$ is a completion of the  metric space $(D(\Rl),\Wrk)$.

\begin{corollary} \label{cor:Lambda-bar-lipschitz}
    $\Lambda$ extends to a $1$-Lipschitz map $\overline{\Lambda}: (\overline{D}(\Rl),\Wrk) \to L^1(\N \times \R)$.
\end{corollary} 



\section{Stability of barcodes} \label{sec:chain-complexes}

Let $C_\bullet$ be a finite dimensional chain complex of $\mathbf{k}$-vector spaces with a specified  basis. 
That is, we have
a finite graded set $S = \bigsqcup_{j \in \Z} S_j$
and $C_\bullet = (\mathbf{k}S_j,d_j)_{j \in \Z}$,
where each $d_j: \mathbf{k}S_j \to \mathbf{k}S_{j-1}$ is a $\mathbf{k}$-linear map and for all $j$, $d_{j-1} \circ d_j = 0$.
For example, $C_\bullet$ may be the simplicial chain complex of a finite simplicial complex or the cellular chain complex of a finite cell complex.
In this section, we will assume that all of our chain complexes are of this form.

Define a relation on $S$ given by $\sigma \prec \tau$ if 
there exists a $j$ such that $\tau \in S_j$, $\sigma \in S_{j-1}$, and the coefficient of $\sigma$ in $d_j(\tau)$ is nonzero.
The reflexive and transitive closure of $\prec$ is a partial order on $S$ which we denote $\subset$.
If $C_\bullet$ is the simplicial chain complex of a simplicial set $S$, then $(S,\subset)$ is the face poset.

Fix a compact interval $[a,b]$.
Assume that we have a function $w: S \to [a,b]$ such that $w(\sigma) \leq w(\tau)$ for all $\sigma \subset \tau$.
We call such a function  a \emph{weight} on $C_\bullet$.
Let $(C_\bullet,w)$ denote the functor from $(\R,\leq)$ to 
the category of sub-chain complexes of $C_\bullet$ given by
$w^{-1}(-\infty,-]$.
For $j \in \Z$, $(C_j,w)$ is a persistence module, and 
\begin{equation} \label{eq:Cj-infty}
    (C_j,w) \isom \bigoplus_{\sigma \in S_j} [w(\sigma),\infty).
\end{equation}


Since the possible values of $w$ are restricted to the codomain $[a,b]$, we also restrict the indexing category of our persistence module to $[a,c)$, where $b < c < \infty$.
In fact, our results will be stronger for smaller values of $c$, and we will determine the minimum $c$ for which they hold.
From now on, our persistence modules are assumed to be indexed by $[a,c)$ with $b < c < \infty$.
Under this restriction, \eqref{eq:Cj-infty} becomes
\begin{equation} \label{eq:Cj-c}
    (C_j,w) \isom \bigoplus_{\sigma \in S_j} [w(\sigma),c).
\end{equation}
Furthermore, for $j \in \Z$, we have the persistence module $H_j(C_\bullet,w)$ given by composing $(C_\bullet,w)$ and the functor $H_j$, given by taking homology in degree $j$.
Let $B_j(C_\bullet,w)$ denote the barcode of $H_j(C_\bullet,w)$.
Let $B(C_\bullet,w)$ denote the graded set $\{B_j(C_\bullet,w)\}_{j \in \Z}$.

Given weights $w$ and $v$ on $C_\bullet$ with values in $[a,b]$, we define 
\begin{equation} \label{def:drank-chain-group}
    \drank((C_j,w),(C_j,v) = \sum_{\sigma \in S_j} \drank([w(\sigma),c),[v(\sigma),c), \text{ and}
\end{equation}
\begin{equation*}
    \drank((C_\bullet,w),(C_\bullet,v)) = \sum_j \drank((C_j,w),(C_j,v)).
\end{equation*}
We also define
\begin{equation*}
    \Wrk(B(C_\bullet,w),B(C_\bullet,v)) = \sum_j \Wrk(B_j(C_\bullet,w),B_j(C_\bullet,v)).
\end{equation*}

In this section we prove the following main results.
\begin{theorem} \label{thm:stability-chains}
    For all chain complexes $C_\bullet$ and all weights $w$ and $v$ with values in $[a,b]$ and corresponding persistence modules indexed by $[a,c)$ where $c > b$, 
    \[
        \Wrk(B(C_\bullet,w),B(C_\bullet,v)) \leq \drank((C_\bullet,w),(C_\bullet,v)).
    \]
    if and only if $c - b \geq b-a$.
\end{theorem}

\begin{theorem} \label{thm:stability-j-chains}
    For all chain complexes $C_\bullet$ and all weights $w$ and $v$ with values in $[a,b]$ and corresponding persistence modules indexed by $[a,c)$ where $c > b$, 
    \[
        \Wrk(B(C_j,w),B(C_j,v)) \leq \drank((C_j,w),(C_j,v)) + \drank((C_{j+1},w),(C_{j+1},v)).
    \]
    if and only if $c -b \geq b-a$.
\end{theorem}

These two results are special cases of the following.
\begin{theorem} \label{thm:stability-barcode}
    Let $m \leq n$.
    For all chain complexes $C_\bullet$ and all weights $w$ and $v$ with values in $[a,b]$ and corresponding persistence modules indexed by $[a,c)$ where $c > b$, 
    \[
        \sum_{j=m}^n \Wrk(B_j(C_\bullet,w),B_j(C_\bullet,v)) \leq \sum_{j=m}^{n+1} \drank((C_j,w),(C_j,v))
    \]
    if and only if $c \geq 2b-a$.
\end{theorem}

We observe that the right hand side of the inequalities above are smaller if $c$ is smaller, so the results are strongest for $c= b + (b-a) = 2b-a$.

Before proving \cref{thm:stability-barcode}, we prove the following preliminary results.

\begin{lemma} \label{lem:path-weights}
    Let $w$ and $v$ be weights on $C_\bullet$.
    Then for $0 \leq t \leq 1$, $(1-t)w + t v$ is a weight on $C_\bullet$.
\end{lemma}

\begin{proof}
    This follows from the definition of a weight on $C_\bullet$.
\end{proof}

\begin{lemma} \label{lem:path-weights-geodesic}
    Let $c,v_0,v_1 \in \R$ with $v_0,v_1 < c$.
    Consider the intervals $[v_0,c)$ and $[v_1,c)$.
    Then the path of intervals $\gamma_t = [v_t,c)$ is a $\drank$-geodesic from $[v_0,c)$ to $[v_1,c)$.
\end{lemma}

\begin{proof}
    It is easy to check that for $0 \leq t \leq 1$, $\drank([v_0,c),[v_t,c)) + \drank([v_t,c),[v_1,c)) = \drank([v_0,c),[v_1,c))$.
\end{proof}

\begin{proposition} \label{cor:geodesic-pm}
    Let $w_0$ and $w_1$ be weights on $C_\bullet$ with values in $[a,b]$.
    Let $w_t = (1-t)w_0 + t w_1$.
    Consider persistence modules to be indexed by $[a,c)$ where $c > b$.
    Then for $j \geq 0$,
    the path $(C_j,w_t)$ of persistence modules is a $\drank$-geodesic from $(C_j,w_0)$ to $(C_j,w_1)$.
\end{proposition} 

\begin{proof}
    By \cref{lem:path-weights}, $w_t$ is a path of weights on $C_\bullet$.
    Thus, for $j \geq 0$,  we have a one parameter family of persistence modules $(C_j,w_t)$.
    By \eqref{def:drank-chain-group} and \cref{lem:path-weights-geodesic},
    $(C_j,w_t)$ is a $\drank$-geodesic.
\end{proof}

\begin{lemma} \label{lem:ineq}
    $\drank([x,y),[x',y')) \leq \frac{1}{2}\abs{x-x'}(b-a) + \frac{1}{2}\abs{y-y'}(b-a)$.
\end{lemma}

\begin{proof}
    There are three cases to consider: $x < y \leq x' < y$, $x \leq x' < y' \leq y$, and $x<x'<y<y'$. 
    In each of these the cases, the support of $\abs{\rank [x,y) - \rank [x',y')}$ is contained in $[x,x'] \times [a,b] \cup [a,b] \times [y,y']$.
    The result follows from \cref{eq:drank}.
\end{proof}

\begin{lemma} \label{lem:iff}
    If $c \geq 2b-a$ then
    for all $x,x' \in [a,b]$, $\frac{1}{2} \abs{x-x'}(b-a) \leq \drank([x,c),[x',c))$.
\end{lemma}

\begin{proof}
    Assume that $c \geq 2b-a$.
    Then $c-b \geq b-a$.
    Assume $x < x'$.
    The support of $\abs{\rank [x,c) - \rank [x',c)}$ contains $[x,x') \times [b,c)$.
    Thus $\drank([x,c),[x',c)) \geq \frac{1}{2} \abs{x-x'}(b-a)$.
\end{proof}


\begin{proposition} \label{cor:drank}
    If $c \geq 2b-a$ then
    for all $[x,y], [x',y'] \subset [a,b]$,
    \[
        \drank([x,y),[x',y')) \leq \drank([x,c),[x',c)) + \drank([y,c),[y',c)).
    \]
\end{proposition} 

\begin{proof}
    Assume $c \geq 2b-a$.
    Then the desired inequality follows from \cref{lem:ineq} and \cref{lem:iff}.
%
\end{proof}

\begin{proof}[Proof of \cref{thm:stability-barcode}]
    Assume $c \geq 2b-a$.
    Consider weights $w_0$ and $w_1$ on $C_\bullet$. Let $j \geq 0$.
    By \cref{lem:path-weights}, $w_t = (1-t) w_0 + t w_1$ is a path of weights on $C_\bullet$.
    By \cref{cor:geodesic-pm}, $(C_j,w_t)$ is a $\drank$-geodesic.

    


    Since $w_t$ is a PL curve of weights, $B_j(C_\bullet,w_t)$ is a PL vineyard (see the original reference \cite{csem:vineyards} or \cite[Theorem 3.6]{Bubenik:2023a} for details in terms of weights).
    Thus, there is an $n \geq 1$ and $0 = t_0 < t_1 < \cdots < t_n = 1$, such that for each $1 \leq i \leq n$ and for all $t$ in the closed interval $[t_{i-1},t_i]$,  
    $B_j(C_\bullet,w_t) = \{[w_t(\sigma_k),w_t(\tau_k))\}_{k=1}^p \cup \{[w_t(\rho_\ell),c)\}_{\ell=1}^q$, for $j$-cells $\{\sigma_k\}_{k=1}^p \cup \{\rho_\ell\}_{\ell=1}^q$, called positive $j$-cells, and for $(j+1)$-cells $\{\tau_k\}_{k=1}^p$, called negative $(j+1)$-cells. 

    Since $(C_j,w_t)$ is a $\drank$-geodesic, it suffices to consider $t,t' \in [t_{i-1},t_i]$ for $0 \leq i \leq n$.
    Let $w = w_t$ and $v = w_{t'}$.    
    Then by the definition of Wasserstein distance,
    \begin{multline*}
        \Wrk(B_j(C_\bullet,w),B_j(C_\bullet,v)) \leq \\
        \sum_{k=1}^p \drank([w(\sigma_k),w(\tau_k)), [v(\sigma_k),v(\tau_k))) +
        \sum_{\ell=1}^q \drank([w(\rho_\ell),c), [v(\rho_\ell),c) ).
    \end{multline*}
    By \cref{cor:drank}, we have,
    \begin{equation} \label{eq:Wrk}
        \Wrk(B_j(C_\bullet,w),B_j(C_\bullet,v)) \leq \sum_\sigma \drank([w(\sigma),c), [v(\sigma),c) ),
    \end{equation}    
    where the sum is over all positive $j$-cells and negative $(j+1)$-cells.
    Summing \eqref{eq:Wrk} for $m \leq j \leq n$, we obtain the desired result.

    Assume the inequality holds for all chain complexes $C_\bullet$ and all weights $w$ and $v$.
    Consider the chain complex $C_\bullet$ with graded basis $S = S_m \sqcup S_{m+1}$, where $S_m = \{\sigma\}$ and $S_{m+1} = \{\tau\}$.
    Let $\varepsilon \in (0,b-a)$.
    Let $w$ and $v$ be the weights given by 
    $w(\sigma) = a$, $w(\tau) = b$,
    $v(\sigma) = a$, and $v(\tau) = b - \varepsilon$.
    Then the left hand side of the inequality equals 
    $\Wrk([a,b),[a,b-\varepsilon))$, which equals
    $\frac{\varepsilon}{2}(b-a) - \frac{1}{4}\varepsilon^2$,
    and the right hand side of the inequality equals
    $\drank([b,c),b-\varepsilon,c))$, which equals
    $\frac{\varepsilon}{2}(c-b) + \frac{1}{4}\varepsilon^2$.
    Thus $c \geq 2b-a-\varepsilon$.
    Therefore $c \geq 2b-a$.
\end{proof}

\subsection*{Acknowledgments}

This research was partially supported by the National Science Foundation (NSF) grant DMS-2324353.


\printbibliography

\end{document}